\documentclass[12pt, reqno, twoside, letterpaper]{amsart}

\usepackage{paperstyle}
\usepackage{graphicx}

\usepackage{mathtools}
\usepackage{todonotes}
\usepackage[norefs,nocites]{refcheck}
\usepackage{amsmath,amsthm,amscd,amssymb}

\newcommand{\be}{\begin{equation}}
\newcommand{\ee}{\end{equation}}

\title[The alien in the Riemann zeta function]
{The alien in the Riemann zeta function}
      
\author[W.~Banks]{William D.~Banks}

\date{April Fools Day, 2025}

\begin{document}




\maketitle

\thispagestyle{empty}

\bigskip\begin{center}
\includegraphics[width=5.5in]{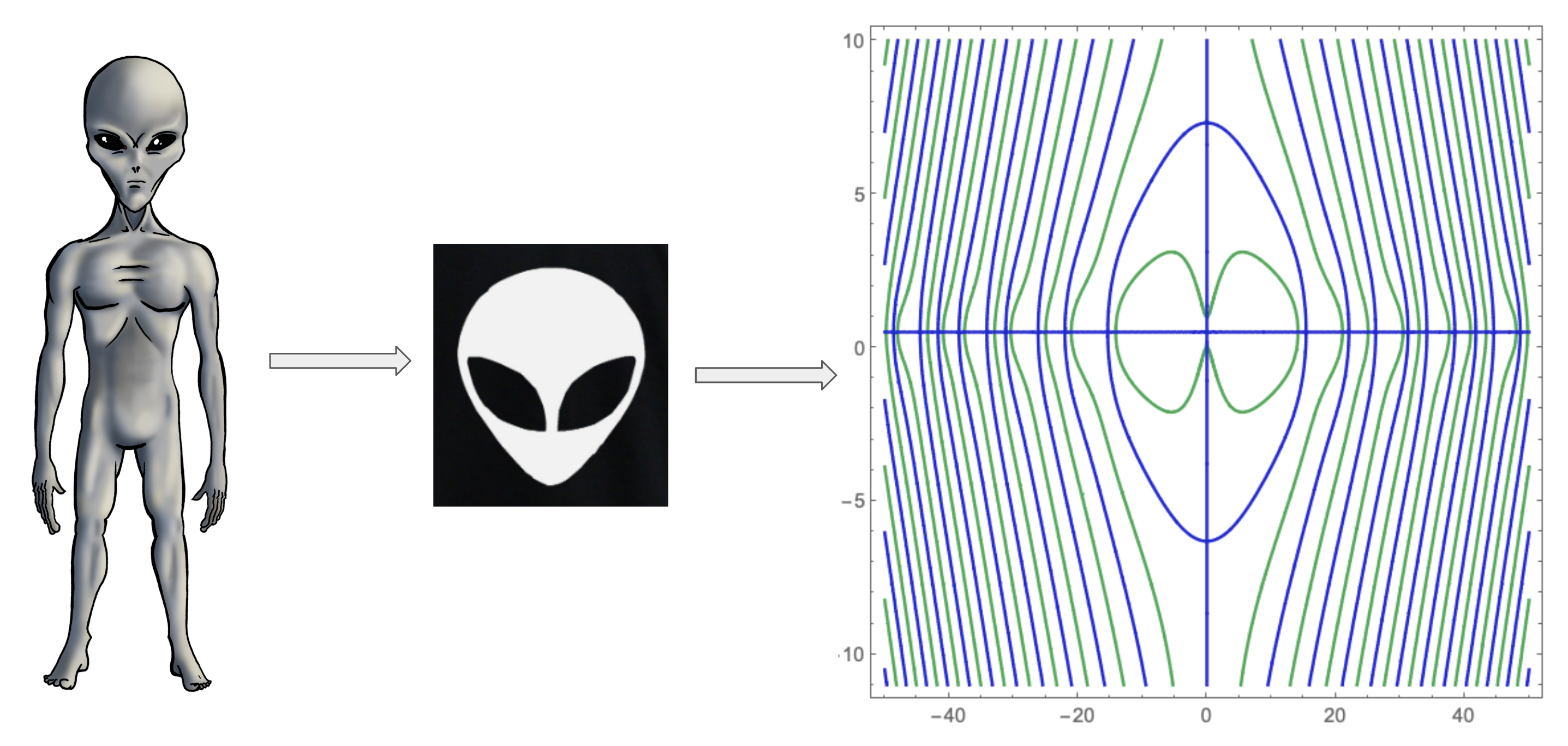}
\end{center}
According to Wikipedia, \textbf{grey aliens}, also known as
\textbf{Zeta Reticulans}, are purported to be extraterrestrial beings; these
 Zetas are frequent subjects of close encounters and alien abduction claims. 
 
The \textbf{zeta function} of Riemann,
\[
\zeta(s)\defeq\sum_{n\in\N}\frac{1}{n^s}=\prod_{p~\text{prime}}
(1-p^{-s})^{-1}\qquad(\sigma>1)
\]
\emph{exhibits its own Zeta Reticulan}.
To produce the rightmost image above, one displays
complex numbers $s=\sigma+it$ with
$\sigma$- and $t$-axes \emph{reversed}. Thus, in the image,
the set of real numbers (i.e., $t=0$) forms the \ccb{\textbf{vertical} blue line},
and the set of numbers $\tfrac12+i\,\R$ (i.e., $\sigma=\tfrac12$, the \emph{critical line}) constitutes the
\ccb{\textbf{horizontal} blue line}.
Set
\[
Z(s)\defeq\pi^{-s/2}\Gamma(s/2)\zeta(s)
\qquad\forall\,s\in\C\setminus\{0,1\}.
\]
and color points in the plane 
\[
\begin{cases}
\ccb{BLUE}&\quad\hbox{if $Z(s)$ is real},\\
\ccg{GREEN}&\quad\hbox{if $Z(s)$ is pure imaginary}.
\end{cases}
\]
Note that the \ccg{green eyes} of the alien are formed from green curves
connecting the poles of $f(s)$ at $0$ and $1$. The eyes are surrounded by 
alien's \ccb{blue head}.
Assuming the \emph{Riemann Hypothesis} and that all zeros of $\zeta(s)$ are simple,
every ``non-eyes'' \ccg{green} curve intersects the \ccb{blue} critical line
in a single point but does not intersect any other \ccb{blue} curve.

\end{document}